\documentclass[twoside,11pt]{article}
\usepackage{amsfonts}
\usepackage{amsmath}
\usepackage{amssymb}
\usepackage[mathscr]{eucal}
 \input epsf

\newtheorem{THEOREM}{Theorem}
\newtheorem{LEMMA}[THEOREM]{Lemma}
\newtheorem{DEFINITION}[THEOREM]{Definition}
\newtheorem{REMARK}[THEOREM]{Remark}
\newtheorem{PROPOSITION}[THEOREM]{Proposition}
\newtheorem{COROLLARY}[THEOREM]{Corollary}

\def\beginabs{\begin{abstract}\noindent}
\def\endabs{\end{abstract}}
\def\beginrem{\begin{REMARK}\rm}
\def\endrem{\end{REMARK}}
\def\beginthm{\begin{THEOREM}}
\def\endthm{\end{THEOREM}}
\def\beginlem{\begin{LEMMA}}
\def\endlem{\end{LEMMA}}
\def\beginprop{\begin{PROPOSITION}}
\def\endprop{\end{PROPOSITION}}
\def\begincor{\begin{COROLLARY}}
\def\endcor{\end{COROLLARY}}
\def\begindef{\begin{DEFINITION}\rm}
\def\enddef{\end{DEFINITION}}

\newcommand{\cab}{,\allowbreak}

\newcommand{\imply}{\relax\allowbreak\ifmmode\;{\Rightarrow}\allowbreak
\;\else\ifhmode\unskip\enskip\fi$\Rightarrow$\allowbreak\enskip\fi}

\newcommand{\bigspace}{\par\ifdim\lastskip<2.5 ex 
\removelastskip\penalty-200 \vskip 2.5 ex plus .3ex minus .3ex\fi}
\newcommand{\fairspace}{\par\ifdim\lastskip<2 ex 
\removelastskip\penalty-150 \vskip 2 ex plus .25ex minus .25ex\fi}
\newcommand{\paraspace}{\par\ifdim\lastskip<1.25 ex 
\removelastskip\penalty-100 \vskip 1.25 ex plus .25ex minus .25ex\fi}
\newcommand{\smallspace}{\par\ifdim\lastskip<1 ex 
\removelastskip\penalty-75 \vskip 1 ex plus .2ex minus .2ex\fi}

\newcommand{\proofspace}{\par\ifdim\lastskip<2 ex
\removelastskip\penalty-150  \vskip 2 ex plus .25ex minus .25ex\fi}

\def\endproofsymbol{\leavevmode
\hbox{\rm\vrule height 1.4ex  width 1.4ex depth .8ex }}
\def\endproofsymbol{\mbox{$\square$}}

\def\endproof{\Tag{\endproofsymbol}\proofspace}
\newcommand{\proof}{\par\noindent{\bf Proof.}\kern.5em}

\newcommand{\rs}{root system}

\newcommand{\ip}{innerproduct}
\newcommand{\ld}{linearly dependent}
\newcommand{\li}{linearly independent}
\newcommand{\ev}{eigenvalue}
\newcommand{\lav}{largest eigenvalue}
\newcommand{\lev}{least eigenvalue}

\makeatletter
\newcommand{\R}{\relax\ifmmode\mathbb R\else
$\m@th\mathbb R$\fi}
\newcommand{\Z}{\relax\ifmmode\mathbb Z\else
$\m@th\mathbb Z$\fi}
\newcommand{\N}{\relax\ifmmode \mathbb N\else
$\m@th\mathbb N$\fi}
\newcommand{\E}{\relax\ifmmode\mathbb E\else
$\m@th\mathbb E$\fi}
\newcommand{\F}{\relax\ifmmode\mathfrak F\else
$\m@th\mathfrak F$\fi}
\newcommand{\fclass}{\relax\ifmmode\mathscr F\else
$\m@th\mathscr F$\fi}
\newcommand{\rclass}{\relax\ifmmode\mathscr R\else
$\m@th\mathscr R$\fi}

\newcommand\lgcover[1]{\relax\ifmmode\left\langle#1\right\rangle
\else$\m@th\left\langle#1\right\rangle$\fi}
\newcommand{\abs}[1]{\relax\ifmmode\left|#1\right|\else
$\m@th\left|#1\right|$\fi}
\newcommand{\cardinal}[1]{\relax\ifmmode{|}#1{|}\else
$\m@th{|}#1{|}$\fi}

\newcommand{\norm}[1]{\relax\ifmmode\left\Vert#1\right\Vert\else
$\m@th\left\Vert#1\right\Vert$\fi}
\def\seq#1#2{\relax\ifmmode#1_1,#1_2\cab\ldots\cab#1_{#2}\else
$\m@th#1_1,#1_2\cab\ldots\cab#1_{#2}$\fi}

\newcommand{\Tag}[1]{\ifvmode\else\unskip\fi
\nobreak\hfil\penalty50 \hskip2em \null
\nobreak\hfil#1\skip@\parfillskip\parfillskip\z@skip
\count@\finalhyphendemerits\finalhyphendemerits\z@\par
\parfillskip\skip@\finalhyphendemerits\count@}

\newdimen\LabeLmargin\LabeLmargin=0pt
\newdimen\LabeLwidth\LabeLwidth=0pt
\newdimen\Hangamount\Hangamount=0pt
\newdimen \ListSpace  \ListSpace=2ex
\newdimen\LabelSpace \LabelSpace=1ex

\newcommand{\setListSpace}{\ifdim\lastskip<\ListSpace\removelastskip
\penalty-100 \vskip\ListSpace plus .15 \ListSpace minus .15 \ListSpace\fi} 

\newcommand{\multilist}[2]{\par\begingroup\setListSpace
\parindent=0pt \LabeLmargin=\Hangamount\setbox0=\hbox{#1}
\advance\LabeLmargin by\wd0 \Hangamount=\LabeLmargin
\setbox0=\hbox{#2}\advance\Hangamount by\wd0 \LabeLwidth=\wd0}

\newcommand{\singlelist}[1]{\par\begingroup\setListSpace
\parindent=0pt \setbox0=\hbox{#1}\LabeLwidth=\wd0 }

\newcommand{\setLabelSpace}{\ifdim\lastskip<\ListSpace
\removelastskip\penalty-75  \vskip \LabelSpace plus .15 \LabelSpace
minus .15 \LabelSpace\fi}

\newcommand{\listitem}[1]{\par\setLabelSpace \hangindent=\LabeLwidth
\hangafter=1 \leavevmode \hbox to\LabeLwidth{#1\hfill}\ignorespaces}

\newcommand{\listend}{\par\setListSpace\endgroup}

\def\@begintheorem#1#2{\trivlist\item[\hskip\labelsep{\bfseries
#1\ #2.}]\itshape} %% This definition is not needed
\def\@seccntformat#1{\csname the#1\endcsname.\quad}
\makeatother

\begin{document}\thispagestyle{empty}
\begin{center}
\large\bf A Method of Classifying All Simply Laced Root Systems\end{center}

 {\tabskip 0pt plus 1fil 
  \halign to \hsize{#\hfil\cr
G. R. Vijayakumar\cr
School of Mathematics\cr
 Tata Institute of Fundamental Research\cr
 Homi Bhabha Road, Colaba\cr
 Mumbai 400\,005\cr
    India\cr\noalign{\vskip .5ex}
 Email: vijay@math.tifr.res.in\cr}}
%% Fax: ++91--22--2280--4610, 4611\quad 
%% Phone: ++91--22--2278--3340\cr

\beginabs A \rs\ in which all roots have
same norm is known as a {\it simply laced root system}.
We present a simple method of  classifying
all simply laced \rs s.
\endabs
\begingroup \ListSpace=1ex \rightskip 0pt plus 1fil
\singlelist{{\bf Keywords:}\enspace}
\listitem{\bf Keywords:}  simply laced root system, base,
euclidean space.\listend\endgroup


\noindent This note is motivated by Chapter 12 of \cite{g-r}, a study
of the class of all graphs with least eigenvalue~$\geqslant-2$:
we present a much simpler and shorter method of deriving Theorem
12.7.4 of \cite{g-r}.
 Let $\N$, $\Z$ and $\R$ be respectively the set of all
positive integers, the set of all integers and the set of all reals.
Let $\E$ be the {\it euclidean space\/} of countably infinite
dimension; i.e., $\E$ is the usual \ip\ space defined on
$\{(r_1\cab r_2\cab \ldots)\in \R^{\N}:
 \sum _{i=1}^\infty r_i^2<\infty\}$;
for any  $x=(r_1\cab r_2\cab \ldots)$ and 
$y=(s_1\cab s_2\cab \ldots)$ which belong to $\E$, their \ip\
$\sum_{i=1}^\infty r_is_i$ is denoted by $\lgcover{x,y}$. We denote
the zero-vector of any subspace of $\E$ by 0 itself. 
Let $S$ be any subset of $\E$. Then the set 
$\{\alpha_1v_1+\cdots+\alpha_nv_n:\ \mathrm{for\ each}\ i\leqslant n,\
\alpha_i\in \Z\ \mathrm{and}\ v_i\in S\}$
is denoted by $Z(S)$; any element (subset)
in (of)  $Z(S)$ is said to be {\it generated\/} by $S$.
 $\widehat S$ denotes the set 
  $\{v\in Z(S): \norm{v}=\sqrt 2\}$.
  We associate with $S$, a graph denoted by $G[S]$:
its vertex set is $S$; two vertices are joined if their \ip\ is
nonzero. If $G[S]$ is connected, then $S$ is called {\it
indecomposable}; otherwise it is {\it decomposable}. Note that 
$S$ is decomposable if and only if
 it has a proper subset $T$ such that for all
$x\in T$ and for all $y\in S\setminus T$, $\lgcover{x,y}=0$.
If $S$ is \li\ and for all distinct $x,y\in S$,
$\lgcover{x,y}\leqslant 0$, then $S$ is called {\it obtuse}.
\paraspace Our object is to produce a method 
of classifying every non-empty finite set 
$X$ in $\E$ such that for all $x,y\in X$, $\lgcover{x,x}=2$,
$\lgcover{x,y}\in \Z$ and $x-\lgcover{x,y}y\in X$. Such a set $X$ is
known as a {\it simply laced \rs\/} in the literature and any element in
$X$ is called a {\it root\/} of $X$. 
When $X$ is indecomposable (decomposable) it is also called
{\it irreducible} ({\it reducible\/}). Thus note that $X$ is a
(disjoint) union of mutually orthogonal irreducible simply laced \rs s.
Henceforth $\Phi$ denotes a simply laced \rs.
A subset $\Delta$ of $\Phi$ is called a {\it base\/} of $\Phi$ if $\Delta$
is obtuse and generates $\Phi$.  \{Though this definition appears to be
different from that of a base of a general \rs\ (see \cite{hu} or
\cite{ka}), it can be shown that for simply laced \rs s, they are
equivalent.\}

\beginrem\label{rem-basic} For any $x\in \Phi$,  
 $-x=x-2x= x-\lgcover{x,x}x\in \Phi$. Let
$a,b\in \Phi$; note that $\lgcover{a,b}=1\imply a-b\in \Phi$
and $\lgcover{a,b}=-1\imply a+b\in \Phi$; since
 $\abs{\lgcover{a,b}}\leqslant \norm a\norm b$
where equality holds only when one root is a scalar multiple of the other,
it follows that if $a\neq\pm b$ then $\lgcover{a,b}\in
\{-1,0,1\}$. 
\endrem
\noindent The following result classifies all simply laced \rs s.
\beginthm\label{thm-main}
 If $\Omega$ is an irreducible simply laced root system, then there
exists an automorphism $\theta$ of\/ $\E$ such that $\theta(\Omega)\in
\{A_n:n\in \N\}\cup\{D_n:n\in \N\
\mathrm{and}\ n> 3\}\cup\{E_6,E_7,E_8\}$.\endthm
\noindent
 Let $\{e_i:i=1,2,\ldots \}$ be an orthonormal basis for $\E$. 
 Then

\smallspace {\halign to \hsize{\hfil$#$\tabskip 0pt&
${}=#$\hfil\tabskip 0pt plus 1fil\cr
 \mathrm{for\ each}\ n\in \N,\
 A_n&\{\pm (e_i-e_j):1\leqslant i<j\leqslant n+1\}\ \mathrm{and}\cr
\mathrm{for\ each}\ n\geqslant2,\
D_n&\{\pm e_i\pm e_j:1\leqslant i<j\leqslant n\}.\cr
E_8&D_8\cup \Bigl\{\frac12
\sum_{i=1}^8\epsilon_ie_i: \epsilon_i=\pm 1\ \mathrm{for}\
i\leqslant8\ \mathrm{and}\
\prod_{i=1}^8 \epsilon _i=1\Bigr\},\cr
E_7&\{v\in E_8 : \lgcover {v,a}=0\}\ \mathrm{and}\cr
E_6&\{v\in E_7: \lgcover{v,b}=0\}\cr
\multispan2where
$a$ and $b$ are two vectors in $E_8$ such that
$\lgcover{a,b}=1$.\hfil\cr}
\paraspace}

\begin{table}
\centerline{\epsfxsize=.95\hsize
\epsffile[37 63 740 539]{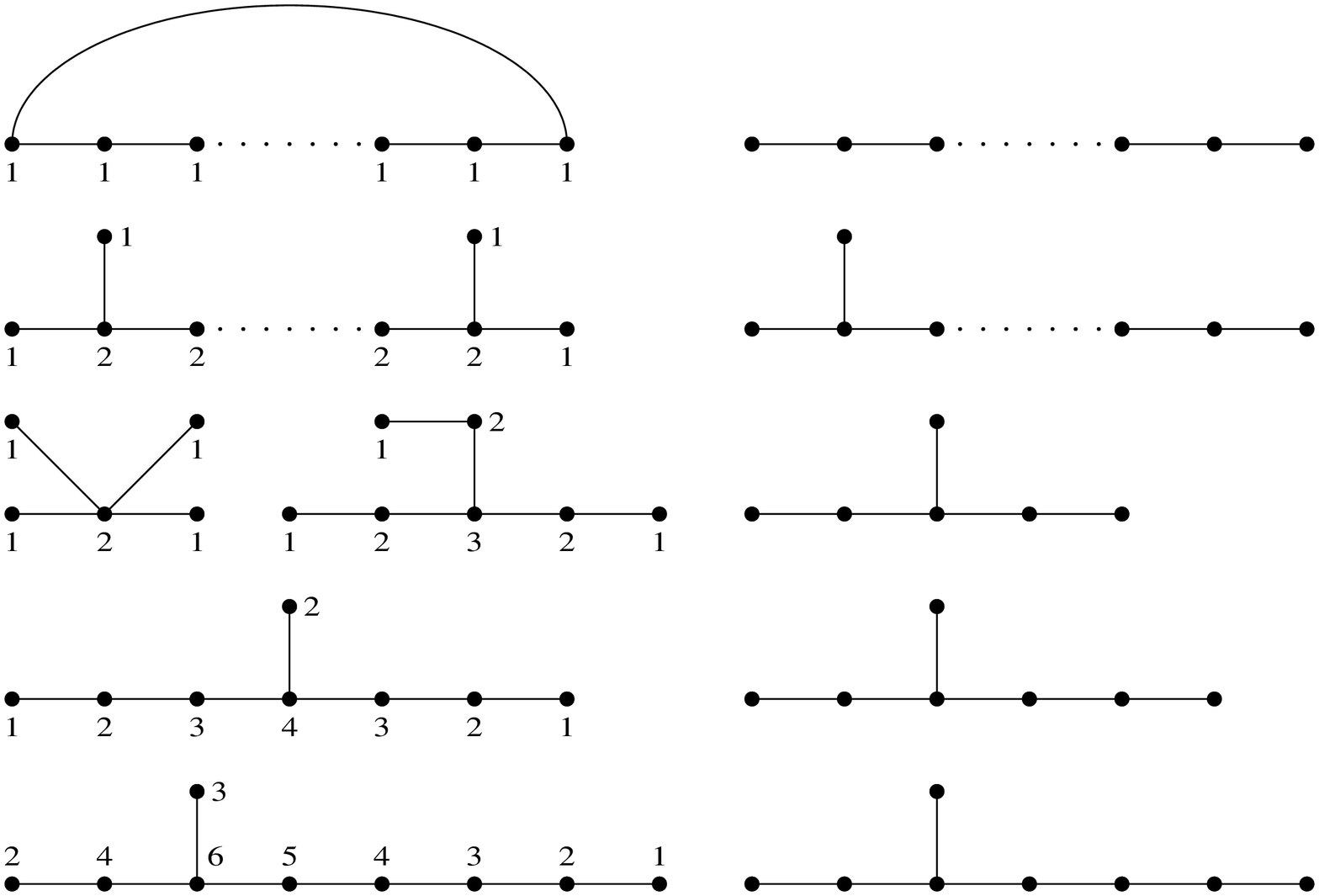}}
\vskip 2ex
\centerline{The family of all connected graphs with \lav\ $\leqslant 2$.}
\end{table}

Let us denote the largest eigenvalue and the least of a graph $G$
by  $\Lambda(G)$ and $\lambda(G)$ respectively.
Derivation of Theorem \ref{thm-main} can be done in two
steps:  (1) For each irreducible \rs, finding a base
$\Delta$ such that $\Lambda(G[\Delta])<2$. 
(2) Associating with each connected graph whose largest 
\ev\ is less than 2, a base of one
of the \rs s defined above. 
The second part is routine;
the reader is 
referred to \cite{hu} and \cite{ka} for this.
We focus only on the first one; our method involves heavily
the family of all connected graphs with largest eigenvalue
$\leqslant2$  displayed in the figure---for each graph $G$ on the
left (right), $\Lambda(G)=2$ ($\Lambda(G)<2$).
(For computation of this family, see \cite{sm}.)

\beginprop \label{prop-base2rs}
 If $X$ is a base of $\Phi$, then $\widehat X=\Phi$.\endprop
\proof For any $w=t_1x_1+\cdots +t_mx_m\in \widehat X$,
where $x_1\cab \ldots \cab x_m$ are distinct vectors in $X$, let
$\rho(w)=\abs{t_1}+\cdots +\abs{t_m}$. It is enough to show
that $\widehat X\subseteq
\Phi$. Let $u=\alpha_1v_1+ \cdots +\alpha_nv_n$ where 
 $v_1\cab \ldots \cab v_n$ are distinct vectors in $X$
 be an element in $\widehat X$; we can assume that
 each $x$ in $\widehat X$ such that $\rho(x)<\rho(u)$ belongs to $\Phi$
and $u\notin X$. Now for some $k\in \{1\cab
\ldots\cab n\}$,  $\alpha_k\ne 0$;  we can assume that
for each $j\leqslant k$, $\alpha_j>0$ and for each $j$ in $\{i\in \N:
k<i\leqslant n\}$, $\alpha_j\leqslant 0$. Then for some $\ell\in \{1\cab
\ldots \cab k\}$, $\lgcover{u,v_\ell}>0$. Since $u\ne v_\ell$, it can
be verified that $\lgcover{u,v_\ell}=1$. Therefore
 $\norm{u-v_\ell}^2=2$; since
$\rho(u-v_\ell)=\rho(u)-1$,
$u-v_\ell\in \Phi$. Now $\lgcover{u-v_\ell,v_\ell}=-1\imply
(u-v_\ell)+v_\ell\in
\Phi$; i.e., $u\in \Phi$.\endproof

\beginprop\label{prop-base}
 $\Phi$ has a base. \endprop
\proof We can assume that $\Phi$ is irreducible.
Let $S$ be an indecomposable
 obtuse subset of $\Phi$ such that $Z(S)\cap\Phi$ is
as large as possible. Suppose that the latter is a proper subset of $\Phi$.
Since $\Phi$ is irreducible, we can find some $p\in \Phi-Z(S)$ and $a\in S$ 
such that $\lgcover{p,a} < 0$.  Let
$T$ be the set of all roots $r$ such that 
 the following holds: for some $x\in S$, $\lgcover{r,x}<0$
and $r$ can be expressed in the form
$p-\sum _{x\in S}\alpha_xx$ where 
for each $x\in S$, $\alpha_x\in \N\cup\{0\}$.
 For any such
$r$, let $\sum_{x\in S}\alpha_x$---it is easy
to verify that this sum is independent of the form for $r$---be
denoted by $\rho(r)$. Note
that $T$ is non-empty because $p\in T$. Choose
a root $q$ in $T$ so that $\rho(q)$ is as large as possible. If can
be easily verified that $\lgcover{q,x}\in\{-1,0\}$ for each $x\in S$.
Now assuming $S\cup \{q\}\subset \R^n$ where $n=\abs{S}+1$, let
$B$ be the $n\times n$ matrix whose rows are the vectors in
$S\cup\{q\}$.  Since $S\cup\{q\}$
is not obtuse, it is \ld. 
Therefore  $B$ is singular. Then the relation
$\lambda(BB^\top)\geqslant 0$ where  $B^\top$ is the transpose
of $B$ becomes an equality. Since
 $BB^\top=2I-A$  where $A$ is the adjacency matrix of $G[S\cup\{q\}]$,
 it follows that $\lambda(2I-A)=0$; therefore
$\Lambda(A)=2$.  Then $G[S\cup\{q\}]$ is one of the graphs on the
left side of the figure.
 It can be
 verified that $\alpha A=2\alpha$ where $\alpha$ is the vector formed by the
labels assigned to the vertices of $G[S\cup\{q\}]$.
Now $(\alpha B)(\alpha B)^\top=\alpha BB^\top
\alpha^\top=\alpha(2I-A)\alpha^\top= 0\alpha^\top =0$. Therefore $\alpha
B=0$.  Since one of the coordinates of $\alpha$ is
1, a vector in $S\cup\{q\}$, say $u$, belongs to $Z(X)$ where
$X=(S\cup \{q\})\setminus \{u\}$. Since $u$ is a pendent vertex of
$G[S\cup\{q\}]$, $G[X]$ is connected; therefore $X$ is indecomposable.
Since $S\subset Z(X)$, $X$ is \li. Note also that $Z(S)\subsetneq Z(X)$ because
$S\cup\{p\}\subset Z(X)$.
Thus the presence of $X$ contradicts the choice of $S$.
 Therefore $S$ is a base of $\Phi$. \endproof

\noindent
Now we can prove the main result: By Proposition \ref{prop-base}, $\Omega$ 
has a base $X$.
Taking $X\subset \R^{\cardinal X}$, let $B$ be the $\cardinal X\times
\cardinal X$ matrix  
whose rows are the elements of $X$. Then $BB^T=(2I-A)$ where
$A$ is the adjacency matrix of $G[X]$. Since $B$ is non-singular,
so is $BB^\top$; therefore $\lambda(BB^\top)>0$; i.e.,
$\lambda(2I-A)>0$. Hence $\Lambda(A)<2$. Since $\Omega$ is
irreducible, $X$ is indecomposable;
therefore $G[X]$ is one of the graphs on the right side of the figure.
Now as mentioned in the discussion before Proposition
\ref{prop-base2rs}, there is a
\rs\ $\rclass\in\{A_n:n\in \N\}
\cup\{D_n:n\in \N\ \mathrm{and}\ n\geqslant 4\}\cup\{E_6,E_7,E_8\}$ having a
base $Y$ such that $G[Y]=G[X]$. Therefore there is a bijection
$f:X\mapsto Y$ such that for all $x,y\in X$,
$\lgcover{f(x),f(y)}=\lgcover{x,y}$. Now the map
$\theta^\star$ from the linear span of $X$
to that of $Y$ defined by 
$\theta^\star\bigl(\sum_{x\in X}\alpha_x
x\bigr)=\sum\alpha_xf(x)$ is an isomorphism. Since these
subspaces are finite dimensional, $\theta^\star$ can be extended to an
automorphism $\theta$ of $\E$. Now by Proposition \ref{prop-base2rs},
$\theta(\Omega)=\theta(\widehat X)=\widehat{\theta(X)}=\widehat Y=\rclass$.
\endproof 
\beginrem\label{rem-rs-generate}
Let $X=\{v_1\cab \ldots\cab v_n\}$ be a subset of $\E$ such that for all 
 $i,j\in \{1\cab \ldots \cab n\}$, $\norm {v_i}=\sqrt 2$ and
$\lgcover{v_i,v_j}\in \Z$.
 For each $x\in \widehat X$, define
$\eta(x)=(\lgcover{x,v_1}\cab \ldots \cab \lgcover{x,v_n})$.
Then $\cardinal{\widehat X}=\cardinal{\{\eta(x):x\in
X\}}$ because for all $a,b\in \widehat X$, $\eta(a)=\eta(b)\imply
a=b$. Therefore  $\widehat X$ is finite. Note for all $a,b\in \widehat X$,
$a-\lgcover{a,b}b\in \widehat X$ because
$\norm{a-\lgcover{a,b}b}^2=\norm{a}^2-2\lgcover{a,b}^2+
\lgcover{a,b}^2\norm{b}^2=2$.
Thus it follows that $\widehat X$ is a simple laced root system.
\endrem 
\noindent Let $A$ be the adjacency matrix of a signed graph  whose
 \lev\ $\geqslant -2$. (By terming
each edge  of a graph as positive or negative, we get a {\it signed
graph}; from the adjacency matrix of the former, that of the latter
is obtained by replacing each entry which corresponds to a negative
edge by $-1$.) Then $\lambda(A+2I)\geqslant 0$. Therefore for some
real matrix $B$, $A+2I=BB^\top$. Thus there is a subset $X=\{v_1\cab
\ldots\cab v_n\}$  of $\E$ such that for all $i,j\in \{1\cab \ldots
\cab n\}$, $\norm {v_i}=\sqrt 2$ and $\lgcover{v_i,v_j}\in \Z$ and
$\left[\lgcover{v_i,v_j}\right]_{i,j=1}^n$---known as the {\it Gram
matrix of $X$}---equals $A+2I$. By Remark \ref{rem-rs-generate}, $X$
is a subset of a simply laced \rs. Therefore by Theorem
\ref{thm-main}, we have the following.

\beginthm If $A$ is the adjacency matrix of a connected signed graph
such that its \lev\ is at least $-2$, then $A+2I$ is the Gram matrix
of a subset of a \rs\ which is  either  $D_n$ for some $n\in \N$
or $E_8$. \endthm

\end{document}